# Information decomposition of symbolic sequences

⊗Eugene V.Korotkov♠♣, Maria A.Korotkova♣ and Nikolai A.Kudryashov♣

♣Moscow Physical Engineering Institute, Kashirskoe Shosse, 31, Moscow, Russian Federation, 115409;
♠Center "Bioengineering" of the Russian Academy of Sciences, 117322, Prospect 60-tya Oktyabrya, 7/1

E.mail:
Eugene Korotkov: katrin22@mtu-net.ru; http://bioinf.narod.ru
Maria Korotkova: katrin2@mail.ru
Nikolai Kudryshov: kudr@dampe.mephi.ru
Fax: 7-095-135-0571
Tel: 7-095-135-2161
⊗Corresponding author

We developed a non-parametric method of Information Decomposition (ID) of a content of any symbolical sequence. The method is based on the calculation of Shannon mutual information between analyzed and artificial symbolical sequences, and allows the revealing of latent periodicity in any symbolical sequence. We show the stability of the ID method in the case of a large number of random letter changes in an analyzed symbolic sequence. We demonstrate the possibilities of the method, analyzing both poems, and DNA and protein sequences. In DNA and protein sequences we show the existence of many DNA and amino acid sequences with different types and lengths of latent periodicity. The possible origin of latent periodicity for different symbolical sequences is discussed.

## 1. Introduction

The development of mathematical methods for the study of symbolical sequence periodicity is given special significance nowadays. Its main concern is the successful determination of DNA sequences from various genomes and the accumulation of a great number of amino acid sequences [1-7]. Therefore, there is a problem for mathematicians and biologists to solve - to determine the structural features of these sequences and to find the biological meaning of the revealed structural features of the sequences.

One such structural feature is a periodicity of symbolic sequences. Earlier comprehensive mathematical methods were developed for the study of periodicity of continuous and discrete numerical sequences, using Fourier transformation and allowing the definition of the spectral density of a numerical sequence [8]. The same approach was widely used for revealing symbolical sequence periodicity. However, such application of Fourier transformation demands presentation of a symbolical sequence as a numerical sequence in which the properties of any symbolical text should be displayed unequivocally. The direct transformation of a symbolical text to a numerical sequence with replacement of symbols by numbers is not an adequate transformation, because we actually introduce weights for symbols, leading to distortion of statistical properties of the initial symbolical sequence. Earlier, several approaches were used for solving this problem [9-17]. The most widely used is the method involving construction from the given symbolical sequence of **m** sequences consisting of the numbers zero and one, formed according to the law: $x(i,j)=1$, if the symbol $a_i$ occupies a site j, and $x(i,j)=0$ in all other cases. Here $\mathbf{A}=\{a_1, a_2, …, a_m\}$ is the alphabet of a symbolical sequence and **m** is the size of the alphabet of a symbolical sequence. The Fourier transformation is then applied to each numerical sequence and the Fourier-harmonics are calculated, corresponding to i-type symbols, as well as matrix structural factors, corresponding to pair correlation of symbols [11]. Final spectral density is usually constructed taking into account statistical characteristics of random sequences [11].

However, in our opinion the given method only works well for the study of periodicity of symbolical sequences with a relatively short length (which is smaller than the size of the symbolical sequence alphabet). For periods with a length greater than the size of the symbolical sequence alphabet, there is the possibility of "decomposition" of the power of the longer periods in favor of the shorter ones. We shall explain this with the following example. Let a symbolical sequence be given with the period YRTDFT repeated 50 times. For this sequence we have 5 numerical sequences consisting of the numbers 0 and 1 (according to the alphabet used). In this case, for the letters Y, R, D, and F the Fourier-harmonics show the length of the symbolical sequence period equal to 6 symbols, but for the letter T the period equal to 3 letters is found. This reduces the power of the 6-letter period by the value of the power of the 3-letter period. This effect will increase with the growth of the relation of period length to the size of the alphabet used. Thus, it turns out that the power of the longer period is a kind of "spread" onto the power of the shorter periods, i.e. there is an effect of attenuation of harmonics with longer periods in favor of harmonics with shorter periods. This effect will be even stronger for cases where there are several replacements in periodic sequences - in such sequences periods could not be simply identical.

The main purpose of the present work is to introduce the concept of information decomposition (ID) of a symbolical sequence that allows the finding of all available cases of periodicity in a sequence, as well as to develop an algorithm which will allow the calculation of ID for a symbolical sequence based on any possible alphabet. In this study, the developed mathematical approach is applied for analysis of poetic texts, DNA sequences and several amino acid sequences. The ID concept allowed us to discover latent periodicity in many poetic texts, genes and various proteins. We also discuss the possible origin of latent periodicity in symbolical sequences of various origin.

## 2. Methods and Algorithms

The main objective is to define the spectral density of a symbolical sequence (or another similar function) without the use of any conversion of a symbolical sequence into a numerical sequence(s). Thus, spectral density of the longer periods should not be distributed onto the shorter periods, as occurs in the case of symbolical sequence decomposition into a certain set of numerical sequences and following the use of Fourier transformation [9-17].

Firstly, let us define the concept of latent periodicity of a symbolical sequence. Let $S = s_1 s_2 \ldots s_k$ be the sequence of symbols with the alphabet $A=\{a_1, a_2, \ldots a_m\}$. If in the sequence $s_1 s_2 \ldots s_k$ there is the subsequence $s_{j+1} s_{j+2} \ldots s_{j+np}$, so that at all *i,r, $1 \leq i \leq n$, $1 \leq r < p$, $s_{j+rp+i} = s_{j+i}$* where *$j \geq 0$, $n \geq 2$, $p \geq 2$,* then the sequence could be presented as $S = s_1 s_2 \ldots s_j [s_{j+1} s_{j+2} \ldots s_{j+n}]_p s_{j+np+1} \ldots s_k$. In this case we can name the fragment $s_{j+1} s_{j+2} \ldots s_{j+n}$ as the period of length **n** repeated in the sequence **p** times. Expanding the concept of periodicity, we name any fragment of the kind *[$s_{j+1} s_{j+2} \ldots s_{j+n}]_p s_{j+1} s_{j+2} \ldots s_{j+l}$,* a periodical fragment, where *l < n*. It is obvious that such periodicity, reflecting some structural features of the symbolical text, is at the same time the simple recurrence of the same fragment. Let us note that the given definition does not describe imperfect periodicity when the replacements of separate symbols in the various periods are possible. Therefore, we shall expand the concept of periodicity further. Let some positions of the period have more than one symbol (set of symbols). Let us introduce the name **"a dim period"** for any period where for each position the allowed set of symbols can be given. Let us name a fragment of a sequence containing, in each appropriate position, one of the allowable symbols of the period, as **a fragment with a dim periodicity**. Let us name as **the latent period** any dim period whose occurrence in the random text is statistically improbable. A set of symbol frequencies for the positions of the period (in view of cyclic rearrangements of positions) we shall name as **a type of a latent period**. For example, the sets of symbols used in various positions of the 6 letter period can look as follows: {YRT} {YR} {DF} {N} {RDF} {YTD}, where the symbols within brackets are symbols used in the corresponding position of the period, and N denotes any symbol. The corresponding set of symbol frequencies or type of latent periodicity is shown in the T matrix:

|   | 1   | 2   | 3   | 4   | 5   | 6   |
|---|-----|-----|-----|-----|-----|-----|
| Y | 1/3 | ½   | 0   | 1/5 | 0   | 1/3 |
| R | 1/3 | ½   | 0   | 1/5 | 1/3 | 0   |
| T | 1/3 | 0   | 0   | 1/5 | 0   | 1/3 |
| D | 0   | 0   | ½   | 1/5 | 1/3 | 1/3 |
| F | 0   | 0   | ½   | 1/5 | 1/3 | 0   |

The elements of this matrix t(i,j) show the probability of the i type symbol in the period position j. This latent period determined by the T matrix of probabilities may be observed in a certain set of symbolical sequences. One such sequence is: [YYDRFT] [RRFDDT] [TYDTFD] [YRDFRY]… From this example, it is obvious that there is no appreciable similarity between the separate periods (placed in square brackets). Revealing such latent periodicity will be rather difficult using Fourier transformation methods because of strong attenuation of harmonics with the longer periods in favor of shorter period harmonics, as was discussed above. Actually, for genetic texts, many latent periods remain unrevealed until now, and the periodicity of many biologically important sequences is not yet revealed. For this purpose it is necessary to develop a mathematical approach allowing the revealing of latent periodicity for periods with rather large lengths, and to define types of latent periodicity, i.e. types of an M matrix. We can expect that there would be various types of matrix M for latent periods with the same length but from different symbolical sequences.

The concept of information decomposition (ID) of a symbolical sequence was applied earlier by authors for analysis of DNA and amino acid sequences [18-22]. Under ID, we understand a spectrum representing the statistical importance of mutual information for periods of various lengths in the analyzed symbolical sequence. Mutual information between the investigated symbolical sequence and artificial symbolical periodic sequences can be used for obtaining an ID spectrum. Let us take a symbolical sequence a(i), i = 1,2, … L. We generate a set of artificial symbolical sequences with periods from 2 up to L/2 of the same length as well as an initial symbolical sequence. We use numbers as symbols of artificial symbolical sequences. The artificial symbolical sequence with period length equal to 2 symbols can be presented as follows: 1,2,1,2,1,2,1,2,1,2,1,2,1,2,…; the sequence with period length equal to 3 symbols can be presented as: 1,2,3,1,2,3,1,2,3,1,2,3,1,2,3,…; the sequence with period length equal to n symbols can be presented as: 1,2,…,n, 1,2,…,n, 1,2…,n…. Further, we can determine the mutual information between the analyzed sequence a(i) and each of the artificial periodic sequences. Values of the mutual information define the ID spectrum for the analyzed symbolical sequence. We fill a matrix M with dimensions (nxk) for the value of the mutual information, where **n** shows the period length of the artificial periodic sequence used, and **k** is the size of the alphabet of the analyzed symbolical sequence. The value of the mutual information is calculated as follows:

$$I = \sum_{1}^{n}\sum_{1}^{k} m(i,j)\ln m(i,j) - \sum_{1}^{n} x(i)\ln x(i) - \sum_{1}^{k} y(j)\ln y(j) + L\ln L \qquad (1)$$

where matrix M shows the numbers of coincidences of ij (i=1,2..,n; j=1,2…k) type between compared sequences (L is the length of the analyzed symbolical sequence, x(i), i=1,2, …, n are the frequencies of symbols 1,2, …, n in the artificial periodic symbolical sequence; y(j), j=1,2, …, k are the frequencies of symbols in the analyzed symbolical sequence). The relation between elements of the T and M matrixes is: t(i,j)x(i)=m(i,j).

One of the properties of mutual information is its orthogonality. This means that I(a, b)=0, in conditions where a and b are symbolical periodic sequences with period lengths representing prime numbers. Another important property of an information spectrum is the nesting of mutual information for various periods one into another [19]. This means that mutual information I(n) for composite period n = $n_1 \times n_2 \times n_3 \ldots \times n_t$ ($n_1$, $n_2$, $n_3$,…,$n_t$ – prime numbers) is equal to:

$$I(n)=I(n_1)+I(n_2)+I(n_3)+\ldots+I(n_t)+I'(n) \qquad (2)$$

where $I(n_1)$, $I(n_2)$, $I(n_3)$,…,$I(n_t)$ are the values of mutual information between artificial periodic sequences with lengths $n_1,n_2,n_3,\ldots,n_t$ and the analyzed symbolical sequence; $I'(n)$ - is "pure" mutual information between the artificial periodic sequence with period length equal to **n** and the analyzed symbolical sequence, which could not arise from periodicity with period lengths equal to $n_1, n_2, n_3, \ldots, n_t$. The formula (2) is easily deduced from the correlation of mutual information for three sequences [23], taking into account that the mutual information between artificial periodic sequences is equal to zero, if the period lengths of them represents various prime numbers.

Formula (2) shows that short ID periods don't provide the effect of "damping" of longer period statistical significance. On the contrary, ID allows accumulation of mutual information of "prime periods" (period lengths are given by various prime numbers) in "composite periods". This ID property is attractive for the detection of long latent periodicity (where the period length value exceeds the value of the alphabet size of the analyzed symbolical sequence).

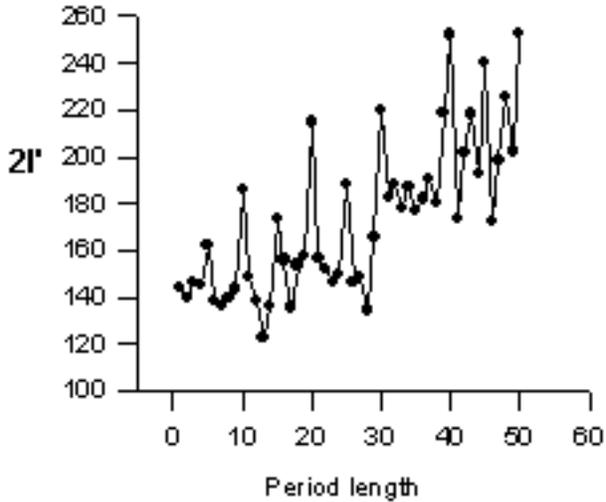

FIG.1. Dependence of the mutual information J(n, 4)=I(n, 4)-3(n-1) at the length of the period for 5-lipoxygenase coding region (34-1000 bases) (Matsumoto et al., 1988) of HUMLOX5A (GENBANK). Mutual information is shown for lengths of the periods multiple to three bases only. On a background of triplet periodicity of DNA is distinctly visible the periodicity in 15 bases that shows herself as local maximums on lengths of the periods equal to 15, 30,45,60,75, … bases.

For ID construction it is necessary to take into account that the value 2I(n,k) is distributed as $\chi^2$ with the value of degree of freedom equal to (n-1)(k-1) [24]. The average value of a mutual information of two random symbolical sequences with alphabet size **n** and **k**, is equal to (n-1)(k-1) correspondingly. This means that a I(n,k) dependence from **n** at the constant **k** has the linear component equal to (k-1)(n-1). Therefore, it is more convenient to show on the diagram the dependence J(n,k)=I(n,k)-(k-1)(n-1), which should resemble a diagram of the Fourier transformation. Such a relation is shown in Fig. 1, and in our previous publications [18,19]. This relation is similar to Fourier harmonics, but only for relatively short periods. For longer periods, the mutual information determined by the formula [1] begins to deviate from the $\chi^2$ distribution if each element of the matrix M becomes less than 10 [25]. This deviation from the $\chi^2$ distribution results in an increase of the average value of J(n,k) with an increasing **n** value. Two approaches can be used for taking into account such a deviation of I(n,k) from the $\chi^2$ distribution in conditions of small sample statistics. The first approach allows direct calculation of the probability of the fact that the relation of symbols in the artificial periodic sequence and in the analyzed symbolical sequence is caused by random factors only, instead of calculation of mutual information [26]. The second approach is based on a Monte-Carlo method for the estimation of the statistical importance of J(n, k) by means of the value Z(n,k) calculation [20-22]:

$$Z(n,k) = \{J(n,k) - \overline{J(n,k)}\} / \sqrt{D(J(n,k))} \tag{3}$$

where $\overline{J(n,k)}$ and D(J(n,k)) show the average value and deviation of the J(n,k) value, for a set of random matrixes with the same sums x(i) and y(j) as in the initial matrix M(n,k). The results of the study of periodicity in various symbolical sequences are presented below. In our study, we used a Monte-Carlo method that permits the execution of the calculations relatively quickly [20-22].

Information decomposition of a symbolical sequence we shall represent as a spectrum $Z(n,k)$. The spectrum $Z(n,k)$ is similar to a spectrum of Fourier transformation for numerical sequences, but has the following advantages: 1. The calculation of the spectrum does not require any transformation of a symbolical sequence to numerical sequences; 2. ID allows the revealing of both the obvious periodicity and the latent periodicity of a symbolical sequence in which there is no statistically important similarity between any two periods; 3. The statistical importance of long periods is not spread onto the statistical importance of shorter periods; 4. On the basis of the matrix M it is possible to determine the type of periodicity.

The proposed approach also allows us to refrain from using a fixed size window when searching for periodicity in a symbolical sequence of large size. For this purpose, we search in the analyzed window for the sub-sequence having the maximum Z value for every tested period length. This approach allows us to realize a segmentation of the analyzed symbolical sequence, depending on the presence of periodicity in various sites of the sequence. Using this approach, periodicity with various period lengths and of various types can be revealed [18-22].

For the study of the long symbolical sequence by the ID method we chose part of a sequence (window) with a length of 2000 symbols. With this window we realized the search for such sub-sequences, where each period length (in a range from 2 up to 200) had the most statistical importance. The window was then moved forward 1000 symbols, and all calculations were repeated. If we analyze the DNA sequence from Genbank or the amino acid sequence from the Swiss-Prot data bank or some linguistic text, it is possible that the length of the analyzed text could be less than the offered window length. In this case, the window size was made equal to the length of the sequence under consideration.

We compared DNA sequences of the ID decomposition method with the Fourier transformation and dynamic programming methods. For calculation of Fourier spectra we used the program "Genescan" [30, 69]. This program converts any DNA sequence to four binary sequences, which can then be Fourier analyzed in the normal manner, to examine correlations between the symbols. For example, the DNA sequence **agatctcggact** converts to four sequences:

```
a    1010000000100
t    0001010000001
g    0100000110000
c    0000101000010
```

The total Fourier spectrum of the DNA sequence is the sum of these individual spectra. We show the power of a periodicity against a period length for the Fourier spectra.

For the search for DNA sequence periodicity with the help of the dynamic programming method we used a program named "Tandem Repeat Finder" [70]. This program used weights for matches, mismatches and indels. These parameters are for Smith-Waterman style local alignment using wraparound dynamic programming. Match weight was +2. Mismatch and indel weights (interpreted as negative numbers) were either 3, 5, or 7. Lower weights allow alignments with more mismatches and indels. The minimum alignment score to report a repeat was chosen as 50.0.

## 3. Results

### 3.1 Model sequence

At first we studied the presence of periodicity in 10 random texts (four-letter alphabet), with the length of each being more than 100 million symbols. It was shown that there were no sub-sequences with Z value exceeding 7.0 in random texts with all lengths of the analyzed periods. Random texts for analysis were generated by the random-number generator of the FORTRAN language using the letter probabilities equal to the probabilities of the bases of the human genome ($p(a)=0.26$; $p(t)=0.26$; $p(c)=0.24$; $p(g)=0.24$). The results obtained demonstrated that we can consider periodicity with Z exceeding 7.0 as statistically important for study of sequences with a usual length of up to 100 million symbols. We used this threshold value for the study of DNA sequences from Genbank. The same result was obtained for alphabets with up to 50 symbols.

The threshold value of Z can be less than 7.0, on the condition that the same level of statistical importance is provided, if we analyze sequences with smaller length. The SWISS-Prot data bank (release 39) contains nearly 39 million amino acids, and the average length of amino acid sequence in the data bank is 392 symbols. It allowed us to choose a threshold value of Z equal to 6.0 for the analysis of amino acid sequences from the SWISS-Prot data bank. The size of poetic texts analyzed did not exceed 30,000 symbols, allowing us to choose a threshold value of Z equal to 5.0 for poetic texts.

We also defined the ID of an artificial symbolical sequence $(ATAAACT)_{100}$. In this sequence there is a distinct period equal to 7 symbols. In Fig. 2 (a) it is shown that the statistical importance (expressed as Z value) is extremely large (exceeding 250) for period length equal to 7 symbols, and gradually decreases for period lengths divisible by 7. So, for such perfect periodicity, the mutual information reaches a maximum value when the period length is equal to 7 symbols. The mutual information maintains this level for period lengths divisible by 7, but the increasing number of degrees of freedom results in the gradual reduction of Z value for the $\chi^2$ distribution.

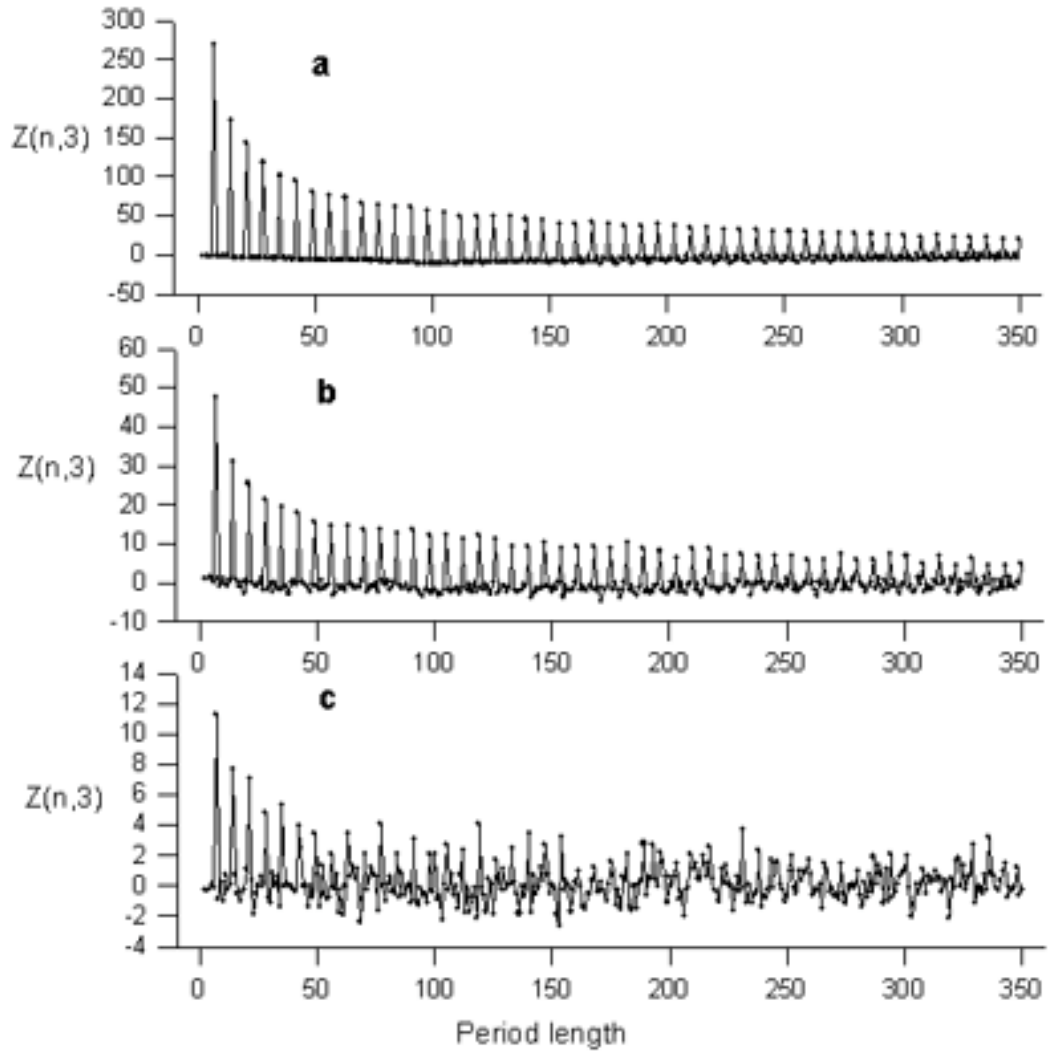

Fig.2. Dependence Z (n,3) from length of the period n for:
**a** - symbolic sequence $(ATAAACT)_{100}$;
**b** - symbolic sequence $(ATAAACT)_{100}$ after introducing of the 50% of random base changes;
**c** - symbolic sequence $(ATAAACT)_{100}$ after introducing of the 80% of random base changes.

We also calculated the ID after the introduction of 50% and 80% random changes in the periodic sequence $(ATAAACT)_{100}$ (Fig. 2 (b,c)). In these figures one can see that the introduction of 50% random changes into the periodic sequence reduces the level of Z(7,3) to 50.0, and the introduction of 80% random changes gives a level of Z(7,3) equal to 11.8. It is possible to determine the period equal to 7 symbols on the statistically important level (Z> 7.0) for the introduction of up to 82% random replacements into the sequence (ATAAACT)$_{100}$. The statistical importance becomes insignificant for periods with length equal to 7 symbols in the case of a greater percentage of introduced letter changes.

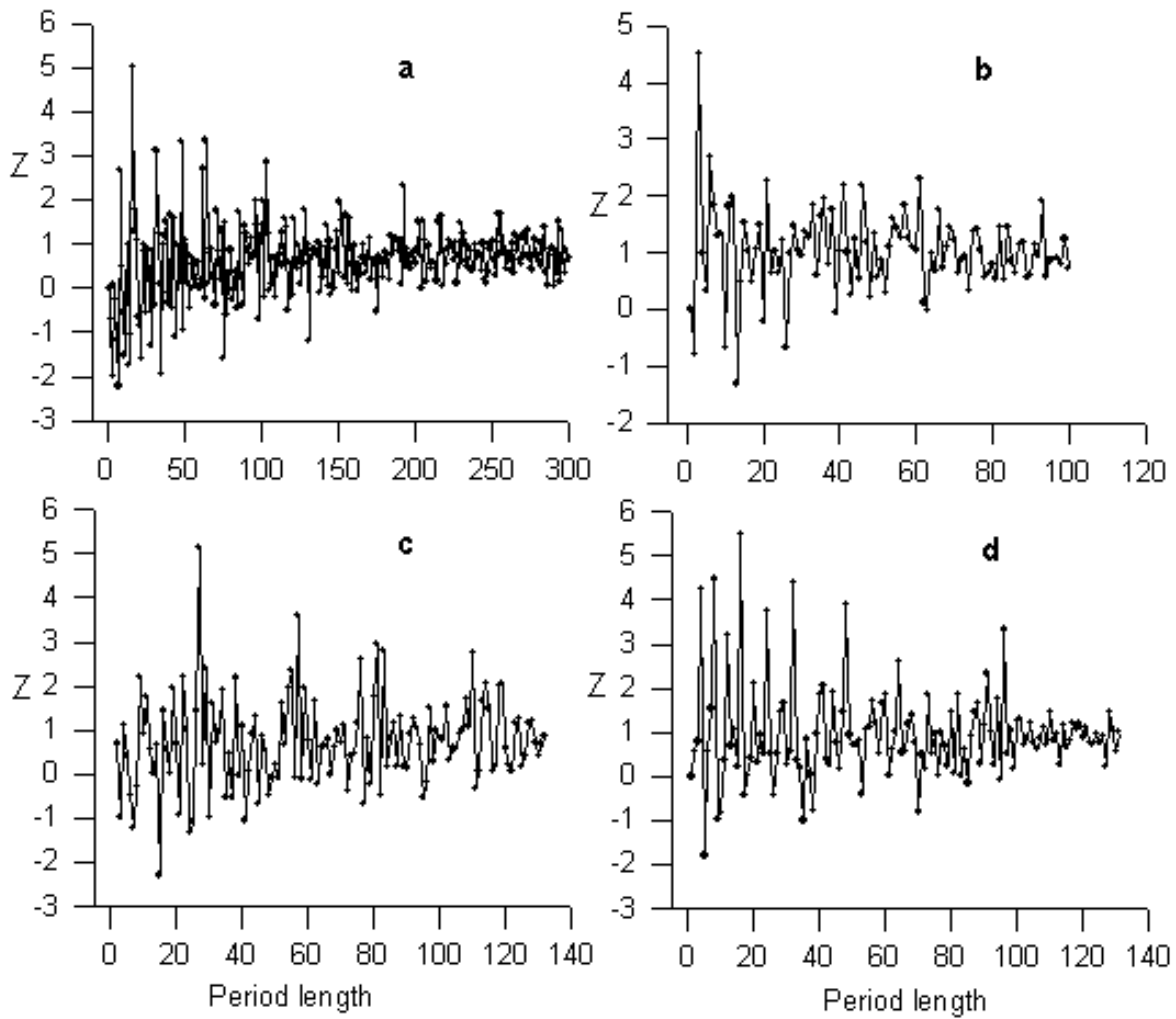

FIG.3. Information decomposition of some poetic texts.
**a** – complete text of the poem of the A.S.Puskin "I remember a wonderful moment… ";
**b** - fragment of the poem V.Mayakovski "Beauties";

    В смокинг вштопорин    Побрит что надо
    По гранд по опере    Гуляю грандом
    Смотрю в антракте    Красавка на красавице
    Размяк характер    Все мне нравиться
    Талии губки ногти    В глянце
    Крашеные губки    Розой убиганятся

**c** - fragments of the products W. Shakespeare "A Midsummer-Night's Dream"
    Have more than thou showest    Speak less than thou knowest
    Lend less than thou owest    Ride more than thou goest
    Learn more than thou trowest    Set less than thou throwest
    Leave thy drink and thy whore    And keep in door
    And thou shalt have more    Thanm two tens to a score

**d** - fragments of the products W. Shakespeare "King Lear"
    If we shadows have offended    Think but this and all is mended
    That you have but slumbered here    While these visions did appear
    And this wear and idle theme    No more yiedling but a dream
    Gentles do not reprehend    If you pardon we will mend
    And as I am an honest Puck

```
                         A                   B                          C
       1      5       10         15        1           1      5       10         15         20         25
  а  1 2 2 2 0 0 0 2 1 2 3 2 0 0 0 1    5 6 7    a  0 1 0 1 1 0 1 1 1 0 0 0 1 2 4 1 0 0 0 0 1 0 0 0 2 0 0
  б  1 1 1 2 2 2 0 0 0 1 0 1 2 0 0 2    0 2 2    b  0 0 0 0 0 0 0 0 0 0 0 0 0 0 0 0 0 0 0 0 0 0 0 0 0 0 0
  в  2 1 0 1 2 0 5 0 3 3 0 3 0 2 2 3    4 4 0    c  0 0 0 0 0 0 0 0 0 0 0 0 0 0 0 1 0 0 0 0 0 0 0 0 0 0 0
  г  2 0 0 0 1 1 1 0 1 0 0 2 1 1 0 1    5 0 3    d  0 0 0 0 0 1 0 1 0 0 1 0 0 0 0 0 1 0 0 1 0 0 1 0 0 0 1
  д  0 3 1 0 0 1 0 1 1 1 2 0 2 1 1 0    1 2 0    e  1 0 1 2 2 3 1 0 4 1 2 0 1 1 0 1 0 0 0 1 1 0 0 0 2 1 2
  е  6 1 7 6 5 3 7 2 0 3 2 1 5 6 2 1    2 7 1    f  0 0 0 0 0 0 0 0 0 0 0 0 0 0 0 0 0 0 0 0 0 0 0 0 0 0 0
  ё  0 0 0 0 0 0 0 0 0 0 0 0 0 0 0 0    0 0 0    g  0 0 0 0 0 0 0 0 0 0 0 0 0 0 0 0 0 0 0 0 0 1 0 0 0 0 0
  ж  1 0 1 2 0 2 0 0 0 2 1 1 0 0 0 0    0 0 0    h  1 0 2 0 0 0 0 1 0 0 0 2 1 4 1 0 1 0 4 1 0 1 0 2 0 1 0
  з  0 2 1 1 1 1 0 2 0 0 0 2 0 0 2 1    0 0 2    i  0 0 0 1 0 0 0 0 0 0 0 0 0 0 0 0 2 0 0 0 0 0 0 0 0 0 0
  и  3 1 2 2 1 3 1 5 0 5 2 2 0 4 5 1    3 5 3    j  0 0 0 0 0 0 0 0 0 0 0 0 0 0 0 0 0 0 0 0 0 0 0 0 0 0 0
  й  0 0 0 2 0 0 2 0 2 0 0 1 0 1 0 3    0 1 0    k  0 0 0 0 0 1 0 0 0 0 0 1 0 0 0 0 0 0 1 0 0 0 1 0 0 0 0
  к  2 2 2 0 0 0 0 1 1 2 0 1 1 0 0 0    4 1 6    l  0 1 0 1 1 0 0 2 1 1 0 0 0 0 0 0 0 0 0 0 0 0 0 0 0 0 0
  л  0 5 1 2 1 2 2 0 2 1 2 0 2 1 1 1    1 0 2    m  0 0 0 0 0 1 0 2 0 0 0 0 0 0 0 1 0 0 0 0 0 0 0 0 0 0 1
  м  1 0 0 0 3 0 1 0 1 1 0 1 1 1 1 4    4 0 1    n  0 0 0 0 0 2 0 1 1 0 0 0 0 1 0 4 1 2 0 0 0 1 0 1 0 2 0
  н  0 3 0 3 5 3 3 1 3 2 0 4 3 2 3 4    7 2 4    o  0 0 0 3 0 0 1 0 2 0 0 1 0 0 0 0 0 2 1 4 2 1 1 2 2 1 0
  о  4 1 4 3 2 7 1 3 2 2 3 1 8 2 4 1    4 8 2    p  0 0 1 0 0 0 0 0 0 0 0 0 0 0 1 0 0 0 0 0 0 0 0 0 0 0 0
  п  1 0 1 0 0 0 1 0 0 3 0 0 1 0 0 0    2 1 2    q  0 0 0 0 0 0 0 0 0 0 0 0 0 0 0 0 0 0 0 0 0 0 0 0 0 0 0
  р  0 1 3 0 0 0 1 2 1 1 2 2 0 2 0 0    7 5 3    r  0 0 1 0 2 0 0 1 0 2 0 0 0 0 1 0 0 1 1 0 0 1 1 1 0 0 0
  с  3 3 0 3 3 0 1 2 0 0 2 3 3 0 0 2    1 0 6    s  2 2 0 1 0 0 1 0 1 2 3 1 0 0 0 1 0 0 0 0 1 0 0 0 0 2 1
  т  1 4 2 1 2 2 3 0 3 1 2 3 0 3 0 0    4 4 3    t  1 4 1 0 0 2 0 0 0 0 4 0 4 1 0 1 0 4 1 0 0 0 3 0 1 0 2
  у  0 0 2 0 0 1 0 0 0 2 2 2 1 0 1 0    1 3 0    u  0 0 0 0 1 0 0 0 0 0 0 0 0 0 0 0 0 1 0 4 1 0 0 0 0 0 0
  ф  0 0 0 0 0 0 0 0 0 0 0 0 0 0 0 0    0 0 0    v  0 0 1 0 0 0 1 0 0 0 0 0 0 1 0 0 0 0 0 0 0 0 0 0 0 0 0
  х  0 0 0 1 1 0 0 0 1 0 1 0 0 0 1 0    0 1 0    w  0 1 1 0 0 0 0 0 0 0 0 0 0 0 0 0 0 0 0 0 0 0 0 1 1 2 1
  ц  0 0 0 0 0 0 0 0 0 0 0 0 0 0 0 1    2 0 0    x  0 0 0 0 0 0 0 0 0 0 0 0 0 0 0 0 0 0 0 0 0 0 0 0 0 0 0
  ч  0 1 0 1 0 1 0 0 2 0 3 0 0 0 0 0    1 0 0    y  0 0 0 0 0 0 0 0 0 0 0 0 0 1 0 0 0 0 0 0 0 0 0 0 0 0 1
  ш  0 0 0 1 0 0 0 0 0 0 1 2 0 0 0 0    1 0 1    z  0 0 0 0 0 0 0 0 0 0 0 0 0 0 0 0 0 0 0 0 0 0 0 0 0 0 0
  щ  0 0 0 0 0 0 0 0 0 0 0 0 0 0 0 0    0 0 0    *  5 1 2 1 2 1 5 1 0 3 1 5 2 1 2 1 5 1 1 2 1 5 1 2 0 0 0
  ъ  0 0 0 0 0 0 0 0 0 0 0 0 0 0 0 0    0 0 0
  ы  0 0 3 0 0 5 0 3 0 1 0 0 3 0 2 0    1 0 0
  ь  2 0 0 1 3 0 2 1 1 1 0 1 2 3 1 1    0 0 0
  э  0 0 0 0 0 0 0 0 0 0 0 0 0 0 0 0    0 0 0
  ю  0 0 0 0 1 0 1 0 0 0 1 0 0 0 0 0    0 2 0
  я  3 0 0 0 0 1 0 3 0 1 3 0 0 1 1 0    3 3 0
  *  5 7 5 4 5 3 61013 3 6 3 3 81111   41018

                         D
       1      5       10         15
  a  0 0 0 0 2 0 1 2 3 0 1 3 0 1 0 5
  b  0 0 1 1 0 0 0 0 0 0 2 0 0 0 0
  c  0 0 0 1 0 0 0 0 0 0 0 0 0 0 0
  d  1 1 1 0 1 1 0 1 1 3 2 1 0 3 1 2
  e  1 1 2 1 3 0 4 4 3 3 3 0 3 1 0 1
  f  0 2 0 0 1 1 0 0 0 0 0 0 0 0 0
  g  0 0 0 0 0 0 1 0 1 0 0 0 0 0 0
  h  0 0 1 0 0 2 1 2 0 1 1 1 0 1 2 1
  i  3 0 2 4 0 0 0 1 1 0 0 1 0 1 1 1
  j  0 0 0 0 0 0 0 0 0 0 0 0 0 0 0
  k  1 0 0 0 1 0 0 0 0 0 0 0 0 0 0
  l  2 3 0 1 1 0 1 0 0 0 0 1 0 0 0 0
  m  0 0 1 0 0 2 2 1 0 0 0 0 1 0 0
  n  0 0 1 0 0 0 0 0 1 4 1 1 0 5 1 1 2
  o  1 0 0 3 1 0 0 0 0 0 1 3 0 0 2 0
  p  0 1 0 0 0 0 0 2 2 0 0 0 0 0 0 0
  q  0 0 0 0 0 0 0 0 0 0 0 0 0 0 0
  r  0 0 1 0 0 1 1 0 0 3 0 2 0 0 0 1
  s  2 0 0 0 2 0 1 1 0 1 0 0 2 1 1 1
  t  1 1 0 0 4 0 1 0 0 0 1 0 1 3 0 1
  u  0 1 1 1 1 1 0 0 0 0 0 2 0 0 0
  v  1 0 0 0 0 0 0 1 0 1 0 0 0 0 0
  w  0 1 0 1 0 0 1 0 0 0 0 1 0 0 2 0
  x  0 0 0 0 0 0 0 0 0 0 0 0 0 0 0
  y  0 0 2 1 0 0 0 0 0 0 0 0 0 0 0
  z  0 0 0 0 0 0 0 0 0 0 0 0 0 0 0
  *  4 6 4 3 0 9 3 0 1 3 5 1 3 3 6 1
```

Fig. 4. Matrixes M for latent periods of the poems that are shown in the Fig3. Letter upper the matrix corresponds to the letter in the Fig.3.

These results show that information decomposition is very stable when a large number of random replacements are introduced. It allows the use of the ID method for the detection of latent periodicity of the symbolical sequence in the absence of any statistically important similarity between any two periods. Such a property of ID is very attractive for the study of DNA sequences, as during evolution DNA sequences accumulated a significant number of base replacements.

### 3.2 Poetic texts.

It seems quite natural to expect that poetic texts will have some degree of periodicity, since such periodicity is quite evident at a perusal of poetic texts. The poetic texts have obvious sound periodicity, which can have a reflection in the structure of the poetic texts. As ID is a sensitive method for the detection of periodicity, it could reveal reflection of sound periodicity of the poetic texts without dependence on the used language. It means that poetic texts in different languages can have various types of periodicity determined by a matrix M, but a common appearance for them should be the presence of any latent periodicity. To notice such generality we analyzed some poetic products in the Russian and English languages.

At first we analyzed several poems of the classic Russian poet, A.S. Pushkin. The ID of the A.S. Pushkin poem "I remember a wonderful moment…" is shown in Fig. 3. In order to realize decomposition, we considered a word space as a letter, and punctuation marks were equated to an additional word space. The maximum value of Z was obtained for a period equal to 16 letters.

So, in matrix M a word space is considered as an additional letter and is marked as *. Matrix M demonstrates that the calculated latent periodicity is conditioned by non-uniform distribution of many letters in period positions (Fig. 4). For example, the Russian letter **E** "prefers" to be located in 1,3,4,5,7 and in 13 positions of the period, the Russian letter **O** is mainly located in 6 and in 13 positions and the Russian letter **B** is mainly located in 7 positions, and so on. The statistical importance of the 16-letter-long period is a result of small deflections from the expected frequency for many letters. Expected values of frequency were calculated for a random text with the same alphabetic structure.

The latent periodicity was also found in other poems of A.S. Pushkin and in poems of other authors. In Fig. 3 (b,c,d) the latent periodicity of fragments in poems of V. Mayakovski and W. Shakespeare is shown. The origin of the periodicity shown in Fig. 3 is also conditioned by preferable distribution of many letters from the Russian or English language alphabets in certain positions of the period. For all four ID shown in Fig. 3 (a,b,c,d) one can see the existence of harmonics of the basic period with the period length divisible by the length of the basic period.

Latent periodicity of poetic texts was revealed in the poems of various authors, writing both in Russian and in English. It is necessary to note that the latent period length in different texts by the same author is not a constant value. Different fragments of text can be characterized by a latent periodicity of various lengths. The type of matrix M (taking into consideration all cyclic rearrangements) for various texts can differ even for periods with identical length.

On the whole, the results obtained from the present study demonstrated that ID is capable of revealing the structure of poetic texts, most probably caused by the sound periodicity of poems in different languages [27,28].

### 3.3 DNA sequences

Earlier, many authors showed the existence of triplet periodicity in DNA sequences in the majority of protein coding regions of genes [11,29-31]. Therefore, triplet periodicity has a coding potential, allowing the prediction of protein coding region positions in many genes [30]. But triplet periodicity can prevent the detection of longer latent periodicity in genes. Since $Z(n,k)$ at **n** lengths equal to 6,9,12,15, … contains $I(3, k)$ as a component of the mutual information, so for the search of periodicity with period lengths divisible by 3, we used random sequences for estimation of $J(n,k)$ and $D(n,k)$ with the same triplet periodicity as in the analyzed DNA sequence. For this purpose we determined frequencies for each of four bases in each of three possible position locations. Then, these frequencies were used for the generation of a set of random sequences with triplet periodicity. Using the random sequences with these properties allowed us to subtract $I(3,k)$ from the mutual information $I(n,k)$ for all of n = 6,9,12,15 …. It is interesting to note that by this method it is possible to take into account the influence of periodicity of any length to longer multiple lengths.

With the algorithm developed, we studied the presence of latent periodicity in various gene sequences from different species from Genbank (release 122). We demonstrated that about 20% of all known genes contain sequences with latent periodicity [18-22]. In Fig. 5 (Ai, Bi, Ci, Di) there are ID of several genes from various bacterial genomes. We chose bacterial genes for demonstration because bacterial genomes contain a low number of the mini- and micro-satellite DNA sequences. It allows us to propose that detected latent periodicity is not connected with recent insertion of different mini- and micro-satellites into genes [32]. Corresponding Fourier spectra are shown in Fig. 5 (Af, Bf, Cf, Df). As it is possible to see from these figures, the Fourier transformation methods do not find the periodicity revealed by the ID method. On the other hand, the program "Tandem Repeats Finder" did not find any periodicity in these DNA sequences.

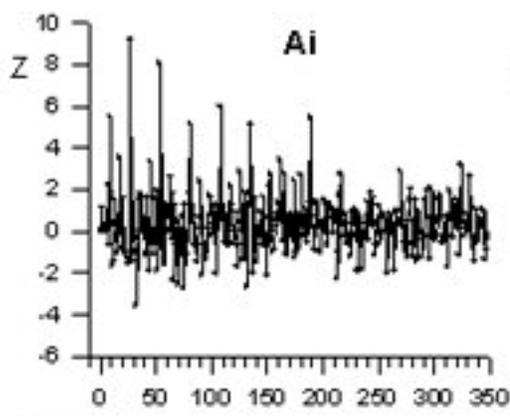
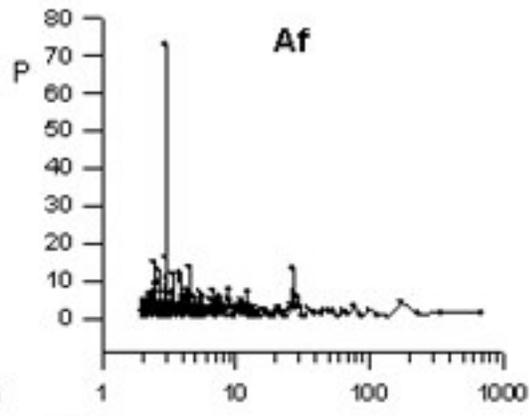
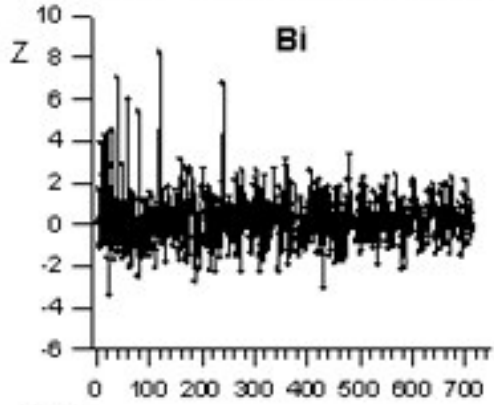
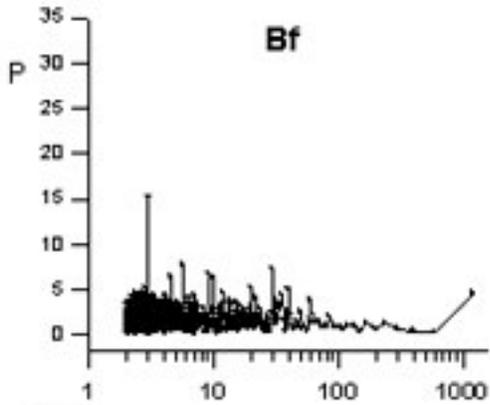
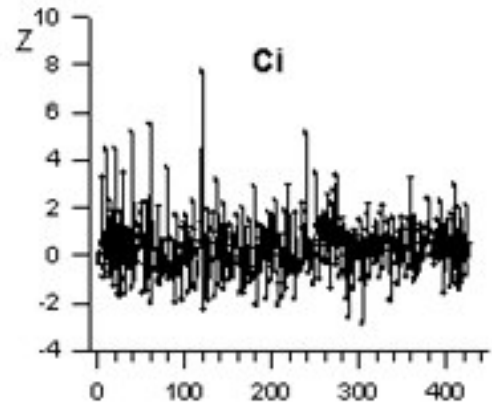
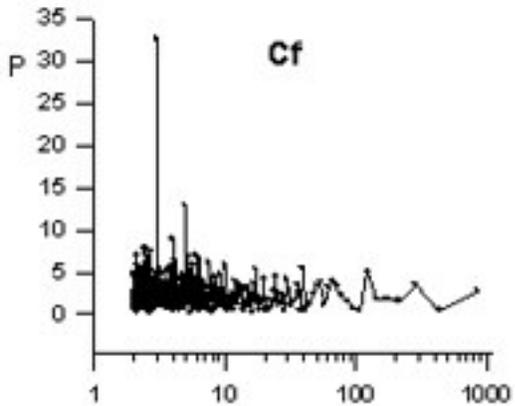
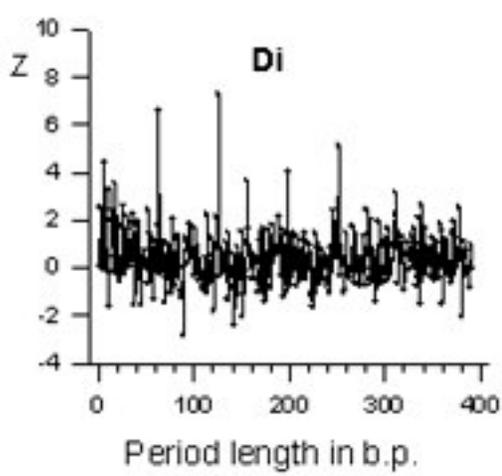
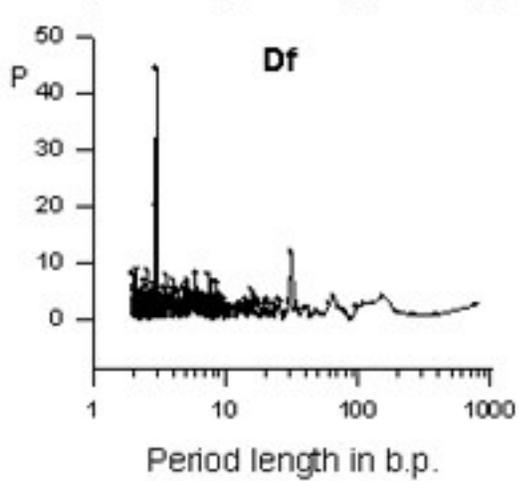

Period length in b.p.    Period length in b.p.

Fig.5. Comparison of the Information Decomposition and Fourier transformation spectrums for the DNA sequences from genomes of different bacteria for period lengths that are multiple to three bases. Z is the Score and P is the Fourier transformation power multiplied by 1000. Figures with letter **i** in the titles show the ID spectrums and figures with letter **f** show the Fourier transformation spectrums. All Fourier transformation spectrums have a strong triplet period and do not show latent periods found by ID methods. It may be caused by the distribution of the power of long latent periods across more short periods by the Fourier method (see main text).

**Ai** –Bordetella pertussisi cya gene for calmodulin-sensitive adenylate cyclase.(981-6101 b.p.) [45] from sequence BPCYA (Genbank). DNA sequence from 5046 to 5742 bases has the latent period equal to 27 bases and Z(27,4) equal to 9.1.

**Bi** - Bacillus subtilis gene [33, 34] for beta-N-acetylglucosaminidase (1296-3938 b.p.) from sequence BACORFX. DNA sequence from 1532 to 2960 bases has the latent period equal to 120 bases and Z(120,4) equal to 9,3.

**Ci** - Deinococcus radiodurans gene [46] for c-di-GMP phosphodiesterase (2867-5239 b.p.) from sequence AE002006. DNA sequence from 3108 to 3963) bases has the latent period equal to 120 bases and Z(120,4) equal to 9,1.

**Di** - Methylobacterium extorquens methanol oxidation gene [47] mxaE (165-1010 b.p.) from sequence AF017434. DNA sequence form 232 to 1015 bases has the latent period equal to 126 bases and Z(126,4) equal to 7,5.

**Af, Bf, Cf** and **Df** - Fourier transformation spectra for the sequences analyzed by ID method in Fig. Ai,Bi,Ci, and Di.

Fig. 5 (Ai) shows the latent periodicity of Bordetella pertussisi cya gene for calmodulin-sensitive adenylate cyclase [45] with period length equal to 27 b.p. It is possible to see from Fig. 5 (Ai, Af) that the latent periodicity equal to 27 bases is observed on a background noise of the strongly expressed triplet periodicity of this gene. The situation reminds us of so-called "amplitude modulation" used in radio engineering for transmission of a signal by a carrier frequency. We found the same situation in many genes from Genbank, and examples of this phenomenon are shown in Fig. 5 (Bi, Ci, Di). Fig. 5 (Ci) shows the latent periodicity of the Deinococcus radiodurans gene for c-di-GMP phosphodiesterase with the length of period equal to 120 b.p. [46], and Fig. 5 (Di) shows the latent periodicity equal to 120 b.p. of a fragment of Methylobacterium extorquens methanol oxidation gene mxaE [47]. We can assume that there are longer correlations in gene structures than triplet correlation and the long correlations can be important for the formation of defined protein structures. Such long correlations can define preferences for using codons in defined gene positions. Triplet periodicity can be considered as some rhythm intrinsic to coding sequences of genes on top of which longer periodicity is expressed.

Fig. 5 (Bi) shows that the major part of the beta-N-acetylglucosaminidase gene [33, 34] has latent periodicity with the period length equal to 120 bases. In Fig. 5 (Bi) one can see the harmonics of the main period with lengths of 240 and 480 bases, as well as induced periodicity for short periods with lengths of 40 and 60 bases. ID can be characterized by the property that the basic periodicity "induces" periodicity on shorter periods, which are constituent elements of the basic period. We consider a period having the maximum Z value in the ID spectrum as the main period. Earlier, a method was developed for estimation of the influence of the main period on the statistical importance of other periods [19].

We also studied the presence of periodicity with period lengths non-divisible by three in DNA sequences from Genbank. We found more than 600,000 such sequences in Genbank. Fig. 6 (Ai, Bi, Ci, Di) shows the examples of revealed latent periodicity and Fig. 6 (Af, Bf, Cf, Df) shows the corresponding Fourier transformation spectra. Latent periodicity found by the ID method can not be revealed in the corresponding Fourier spectra, and any periodicity was not revealed by the program "Tandem Repeats Finder" in these DNA sequences.

The presence of latent periodicity with period length non-divisible by three in gene sequences is very interesting. It can prove that these DNA sequences might have existed long ago as mini-satellite sequences or sequences with multi-tandem duplications, but during evolutionary DNA transformations these DNA sequences might have become a coding sequence. We noticed this phenomenon in the case of a region including two hypothetical proteins of Aeropyrum pernix (Fig. 6 (Bi)), a region of mercuric reductase MerA of Bacillus sp (Fig. 6 (Ci)) and two V-regions from Fugu rubripes of the T-

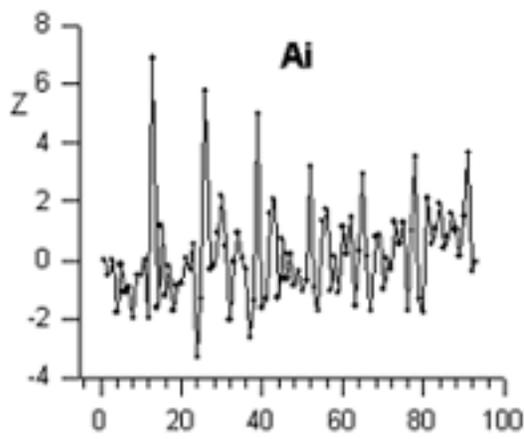
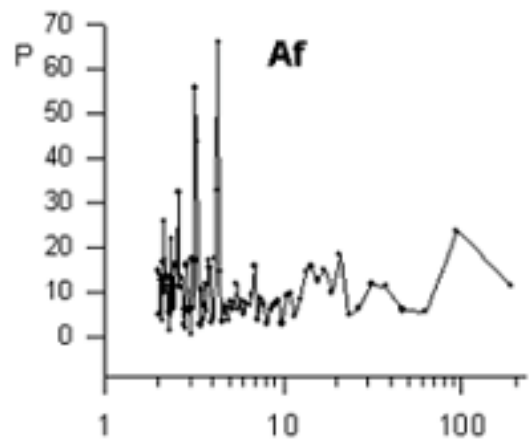
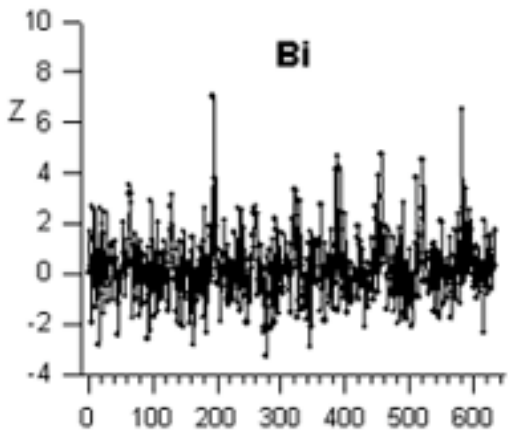
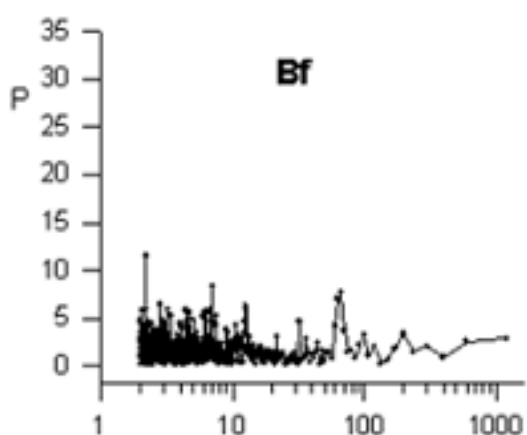
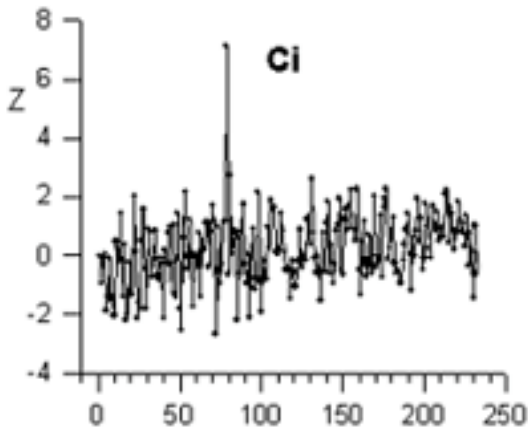
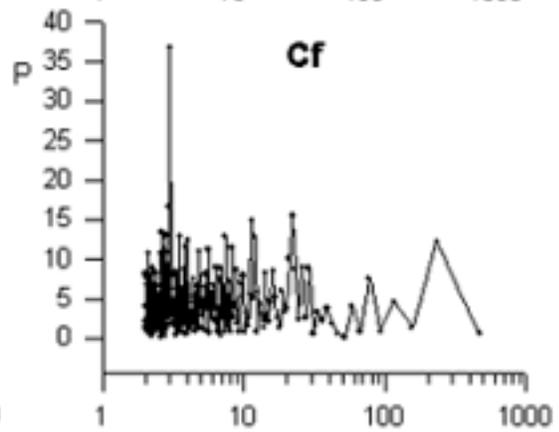
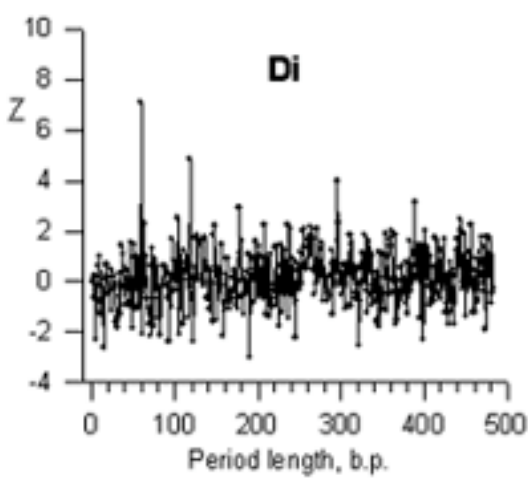
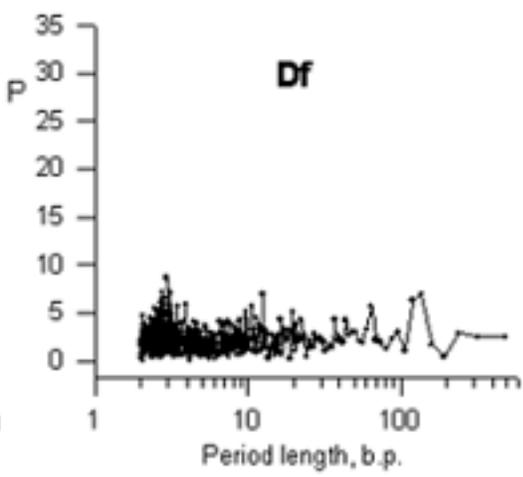

Fig.6. Comparison of the Information Decomposition and Fourier transformation spectrums for the DNA sequences from genomes of different bacteria for period lengths that are unmultiple to three bases. Z is score and P is Fourier transformation power multiplied by 1000. Figures with letter **i** in the titles show the ID spectrums and figures with letter **f** show the Fourier transformation spectrums. Fourier transformation spectra do not show or show very low power for the latent periods found by ID methods. It may be caused by the distribution of the power of long latent periods across more short periods by the Fourier method (see main text).

    **Ai** – noncoding region of the NRC-1 plasmid of Halobacterium sp [48] from AE005163 (Genbank). DNA sequence from 7484 to 7670 bases has the latent period equal to 13 bases and Z(13,4) equal to 7.1.

    **Bi** – region including two hypothetical proteins of Aeropyrum pernix from AP000063 [35]. DNA sequence from 40610 to 41876 bases has the latent period equal to 194 bases and Z(194, base pair 4) equal to 7.2.

    **Ci** – region of mercuric reductase MerA of Bacillus sp. from AF138877 [37]. DNA sequence from 2225 to 2690 bases has the latent period equal to 79 bases and Z(79 ,4) equal to 7.2.

    **Di** – Two V-regions from fugu rubripes of the t-cell receptor alpha-chain (TCRA) gene (13412-13718 b.p., 13850-13892 b.p.) from FRTCRA1 [36]. DNA sequence from 13628 to 14594 bases has the latent period equal to 59 bases and Z(59,4) equal to 7,2.

    **Af, Bf, Cf** and **Df** - Fourier transformation spectra for the sequences analyzed by ID method in Fig. Ai,Bi,Ci, and Di.

cell receptor alpha-chain (TCRA) gene (Fig. 6 (Di)). The DNA sequences from these genes were found with latent period length equal to 59 bases (Fig. 6 (Di)).

In Fig. 6 (Ai) the latent periodicity with length equal to 13 bases, found in non-coding DNA sequence, is also presented. Such period length is usually typical for mini-satellite sequences. Therefore, it is possible to assume that this DNA sequence represents a super-diverged mini-satellite sequence, accumulated after a great many mutations. The ID method, owing to its stability with accumulation of mutations, allowed us to see the periodicity of this sequence. We believe that the ID method can be very useful for the detection of ancient mini- and micro-satellite sequences in various genomes. These sequences can be used for PCR analysis of various genetic traits.

### 3.4 Amino acid sequences.

The latent periodicity of genes should be evident at amino acid level as well, but due to the ambiguity of genetic code several classes of latent periodicity might become statistically insignificant. This means that some types of latent periodicity in genetic sequences could be imperceptible at the level of amino acid sequences, and some types of amino acid latent periodicity could be unnoticeable at the level of DNA sequences. We studied the presence of amino acid sequences with latent periodicity in the Swiss-Prot data bank (release 39). We revealed more than 12,000 amino acid sequences with various types of latent periodicity of various lengths. Figs. 7 and 8 show examples of amino acid sequences with short and long latent periodicity. The periodicity in these examples was observed for protein sequences of different kinds, including various enzymes. This may show that latent periodicity is a characteristic of proteins, fulfilling certain structural functions as well as being catalyzers for certain biochemical reactions. It can be assumed that relatively short latent periods (with period length from 2 up to 10 amino acids) (Fig. 7) are associated with secondary structure formation in amino acid sequences, such as: α-helix, β-layers and some combinations of these structures [32, 38-40]. Longer latent periodicity can be associated with the domain organization of proteins and packing of amino acid sequences [32, 38-40]. In this regard, the long latent periodicity of amino acid sequences can be responsible for hydrogen bond formation between elements of the protein secondary structure, as well as for stabilization of protein structure by means of these interactions.

## 4. Discussion and Conclusions

The developed method of information decomposition (ID) of symbolical sequences proved to be capable of revealing the latent structure of various symbolical sequences. The method intended for detection of periodicity appeared to be rather tolerant to a relatively large number of symbol replacements. The results obtained using this method demonstrated that a large number of known genetic texts contain sequences with latent periodicity of various lengths and various types. It is in agreement with the results of investigation of DNA sequences by integral mathematical methods and finding the long-range correlations (LRC) in DNA sequences [55-68].

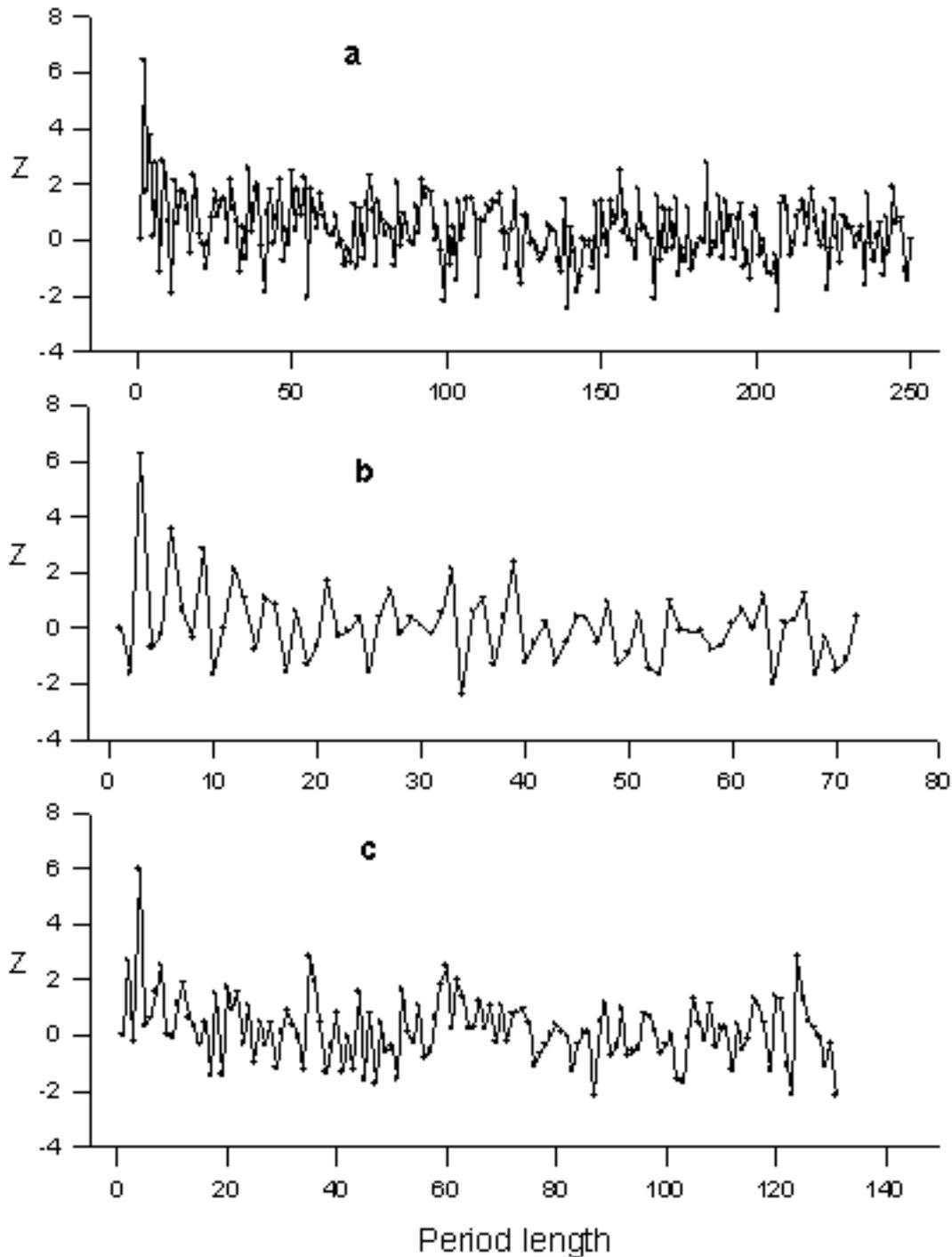

Fig.7. Short latent periodicity of the amino acid sequences.
    **a** - GMP SYNTHASE [49] from Helicobacter pylori (508 amino acids, GUAA_HELPY in Swiss-prot). Sequence from 7 to 505 amino acid has the latent period equal to 2 amino acids and Z(2,20) equal to 6.5
    **b** - ENDOGLUCANASE A PRECURSOR [50] from Butyrivibrio fibrisolvens (429 amino acids, GUNA_BUTFI). Sequence from 64-208 amino acid has the latent period equal to 3 amino acids and Z(3,20) equal to 6.2
    **c** - DNA GYRASE SUBUNIT [51] a Mycobacterium smegmatis (842 amino acids, GYRA_MYCSM). Sequences from 19-280 amino acid has the latent period equal to 4 amino acids and Z(4,20) equal to 6.1

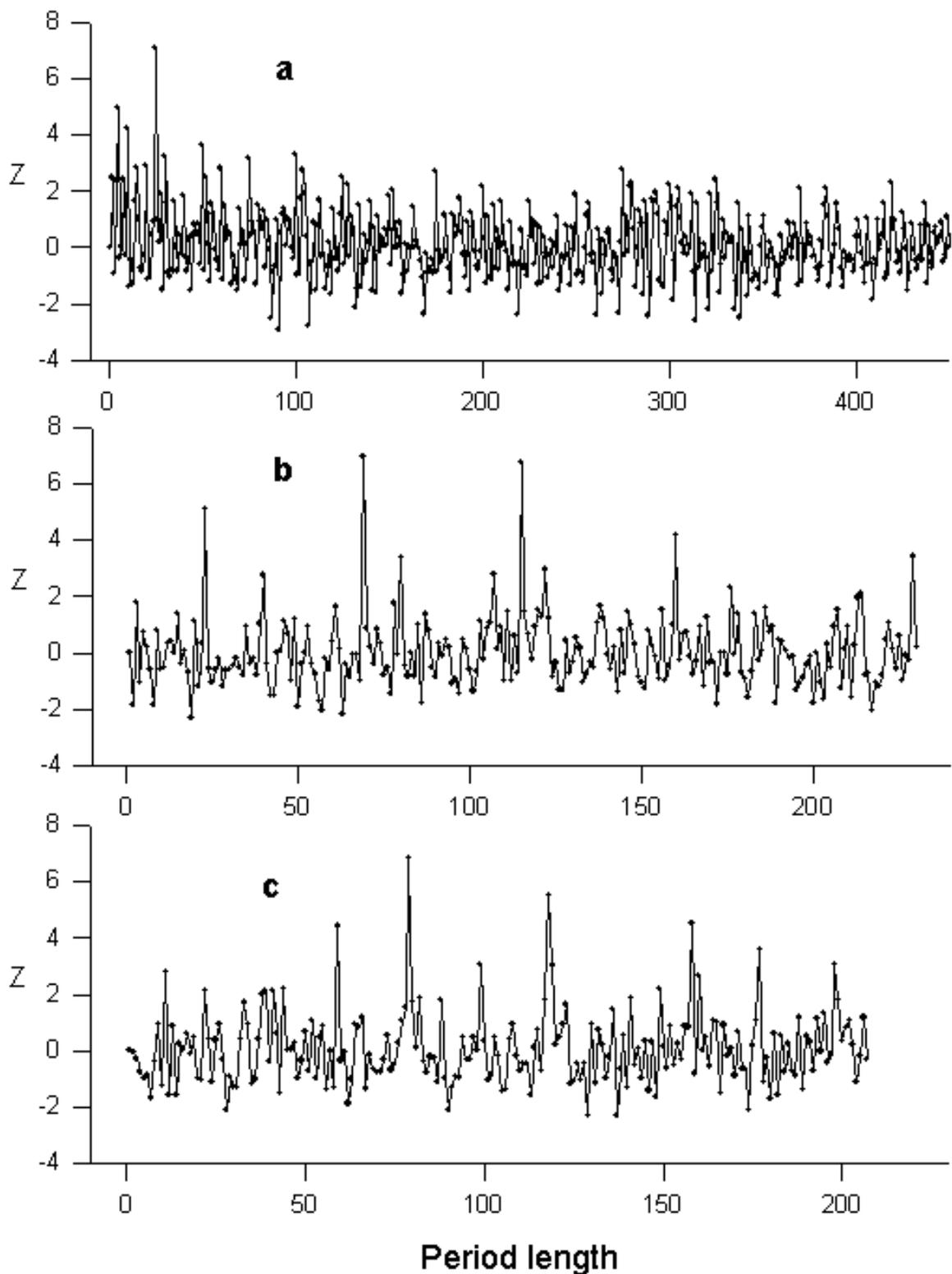

Fig.8. Long latent periodicity of the amino acid sequences.
    **a** - Alpha platelet-derived growth factor receptor precursor [52] from Rattus norvegicus (1088 amino acids, PGDS_RAT in Swiss-prot). Sequence from 85 to 985 amino acids has the latent period equal to 25 amino acids and Z(25,20) equal to 7.1
    **b** - Hypothetical abc transporter ATP-binding protein RV1281C [53] from Mycobacterium tuberculosis (612 amino acids, YAC81_MYCTU) Sequences from 151 to 610 amino acids has the latent period equal to 69 amino acids and Z(69,20) equal to 7.0
    **c** - Putative translational activator C18G6.05C (GCN1 homolog) [54] from Schizosaccharomyces pombe (2670 amino acids, YAC5_SCPO). Sequence from 1421 to 1835 amino acids has the latent period equal to 79 amino acids and Z(79,20) equal to 6.8

It is interesting that a Fourier transformation method [30, 69] can not reveal latent periodicity found by an ID method. The reason is the rather large lengths of found periods, compared with the size of the alphabet used. The power of a long period, as noted above, is distributed on powers of a set of short periods, leading to an impossible detection of latent periodicity. This appearance can be seen by comparing the spectra ID decomposition and Fourier transformation shown in Figs. 5 and 6. It is very interesting also that the method of tandem repeat finding based on dynamic programming [70] was not capable of revealing any periodicity in the sequences shown in Fig. 5 and 6. It can be stipulated that the use of a weight matrix of symbol coincidences limits the possibilities of determinating latent periodicity. Therefore, it seems that the developed programs, based on a dynamic programming method, are only useful for the detection of homologous periodicity. In this sense ID represents a useful method for the definition of rather long latent periodicity in symbolical sequences, which could be missed by other developed mathematical methods and algorithms. We think that the revealed latent periodicity of gene sequences could be caused for two reasons. Firstly, gene sequences could be formed at the time of the first living organisms or some time later by way of simple multi-tandem duplications. After this, natural selection improved their properties over a long time. This improvement could be connected to a set of directed mutations, realized in gene sequences. The directed mutations could change triplet frequencies in genes in correspondence with defined preferences in amino acid usage, that could cause formation of the defined triplet periodicity of gene sequences. However, these mutations probably not could completely destroy the tracks of the processes of gene origin. Therefore, ID shows the presence of a rather long latent periodicity on a background noise of the presence of a triplet periodicity. Secondly, the latent periodicity of the genes which have been found in the present work can be stipulated by defined periodic spatial organization of coded proteins. For example, periodicity equal to 21 bases is usually connected with α-helix formation protein molecules [11, 38, 41, 42]. The longer periodicity could be determined by domain organization formation in proteins [32]; it could also be involved in the process of nucleosome binding with DNA [43]. However, we observed a great variety of period lengths and types in DNA and amino acid sequences. It testifies that some other protein spatial structures have characteristic periods of their own [22]. It is possible to assume that the ID method is able to "see" certain structural characteristics of gene sequences, reflecting spatial organization of the corresponding proteins. In this regard, ID is obviously an important method, allowing the connection of the origin of certain protein structures with the presence of certain latent periodicity in corresponding DNA and amino acid sequences.

As a result of our study, we have found more than $10^6$ sequences with various types of latent triplet periodicity, more than $2 \times 10^6$ DNA sequences with period length from 2 up to 200 bases (without triplet periodicity) in Genbank (release 122) and more than 12,000 cases of latent periodicity of protein sequences in the Swiss-Prot data bank (release 39).

The developed ID method, as well as the methods based on Fourier transformation, are not able to detect latent periodicity with deletions or insertions of symbols. Deletions and insertions of different sizes are mutation events, frequently met in genetic texts. We think that this could lead to a reduction in the number of DNA and protein sequences with latent periodicity found by the ID method in the Genbank and Swiss-Prot data banks. It permits us to assume that the real number of sequences with latent periodicity in genetic texts is actually considerably greater than can be revealed with the ID method at present. This is a reason why we are currently improving the ID method, using some of the methods of profile analysis and dynamic programming. In the near future our results will be completed with a large number of genetic texts with latent periodicity, revealed in the presence of deletions and insertions of symbols.